\documentclass[10pt,twoside,reqno]{amsart}
\usepackage{amssymb}
\usepackage{multirow}
\usepackage{graphicx}
\usepackage[dvips]{epsfig}
\usepackage{latexsym}
\usepackage{amsmath}
\usepackage{amsthm}
\usepackage{amscd}
\usepackage{amssymb}
\input xy
\xyoption{all}

\calclayout

\numberwithin{equation}{section}
\makeatletter

\renewcommand{\@secnumfont}{\bfseries}

\renewcommand{\section}{\@startsection{section}{1}%
  {0mm}{.7\linespacing\@plus\linespacing}{.5\linespacing}
  {\normalfont\bfseries\centering}}

\newcommand{\bibsection}{\@startsection{section}{1}%
  {0mm}{.7\linespacing\@plus\linespacing}{.5\linespacing}
  {\normalfont\scshape\centering}}

\renewcommand{\@biblabel}[1]{#1.}

\theoremstyle{theorem}
\newtheorem{theorem}{\scshape Theorem }[section]

\newtheorem{corollary}[theorem]{\scshape Corollary}
\newtheorem{proposition}[theorem]{\scshape Proposition}
\theoremstyle{definition}

\numberwithin{equation}{section}

\begin{document}

\vspace{1.3cm}

\title {A note on degenerate multi-poly-Genocchi polynomials}

\author{Taekyun Kim$^{1}$}
\address{$^{1}$Department of Mathematics, kwangwoon University, seoul 139-701, Republic of Korea}
\email{tkkim@kw.ac.kr}

\author{Dae San Kim$^{2}$}
\address{$^{2}$Department of Mathematics, Sogang University, seoul 121-742, Republic of Korea}
\email{dskim@sogang.ac.kr}

\author{Han Young kim$^{3}$}
\address{$^{3}$Department of Mathematics, kwangwoon University, seoul 139-701, Republic of Korea}
\email{gksdud213@kw.ac.kr}

\author{Jongkyum kwon$^{*4}$}
\address{$^{4}$Department of Mathematics Education and ERI, Gyeongsang National University, Jinju 52828, Republic of Korea}
\email{mathkjk26@gnu.ac.kr}

\thanks{* is corresponding author}
\maketitle

\begin{abstract}
In this paper, we introduce the degenerate multiple polyexponential functions which are multiple versions of the  degenerate modified polyexponential functions. Then we consider the degenerate multi-poly-Genocchi polynomials which are defined by using those functions and investigate explicit expressions and some properties for those polynomials.
\end{abstract}

\bigskip
\medskip
\section{Introduction}
\medskip

It was Carlitz who initiated the study of degenerate versions of some special numbers and polynomials, namely the degenerate Bernoulli and Euler polynomials and numbers [1]. In recent years, studying degenerate versions of some special polynomials and numbers regained interests of many mathematicians, and quite a few interesting results were discovered [6,7,9,11-15,17].
The polyexponential functions were first introduced by Hardy in [3,4] and rediscovered by Kim [11,14], as inverses to the polylogarithm functions. Recently, the degenerate polyexponential functions, which are degenerate versions of polyexponential functions, were introduced in [11], and some of their properties were investigated. Furthermore, the so-called new type degenerate Bell polynomials were introduced, and some identities connecting these polynomials to the degenerate polyexponential functions were found in [11].
In [12], the (modified) polyexponential functions were used in order to define the degenerate poly-Bernoulli polynomials, and several explicit expressions about those polynomials and some identities involving them were derived.
In [14], the degenerate (modified) polyexponential functions were introduced, and the degenerate type 2 poly-Bernoulli numbers and polynomials were defined by means of those functions. In addition, several explicit expressions and some identities for those numbers and polynomials were deduced. \\
\indent In this paper, we intorduce the degenerate multiple polyexponential functions. These are multiple versions of the degenerate modified polyexponential functions. Then we define the degenerate multi-poly-Gennocchi polynomials by means of those functions. We derive some explicit expressions for the degenerate multi-poly-Gennocchi polynomials and certain properties related to those polynomials.

We recall that, for all $k \in \mathbb{Z}$, the polylogarithm functions are defined by
\begin{equation}\label{01}
\begin{split}
Li_{k}(x) = \sum_{k=1}^{\infty}\frac{x^{k}}{n^{k}}, (\left | x \right | < 1), \quad (\textnormal{see}\,\, [16,18]).
\end{split}
\end{equation}

The polyexponential functions were studied by Hardy in [3,4].\\
Recently, a slightly different version of those functions, which are called the modified polyexponential functions, are defined as an inverse to polylogarithm functions by
\begin{equation}\label{02}
\begin{split}
Ei_{k}(x) = \sum_{n=1}^{\infty}\frac{x^{n}}{(n-1)!n^{k}}, (k \in \mathbb{Z}), \quad (\textnormal{see}\,\,[5,7-9,11,13,14]).
\end{split}
\end{equation}
When $k=1$, by \eqref{02}, we get
\begin{equation}\label{03}
\begin{split}
Ei_{1}(x) = \sum_{n=1}^{\infty} \frac{x^{n}}{n!} = e^{x}-1.
\end{split}
\end{equation}

In [14] (see also [11]), the degenerate modified polyexponential functions are defined by
\begin{equation}\label{04}
\begin{split}
Ei_{k,\lambda}(x) = \sum_{n=1}^{\infty}\frac{(1)_{n,\lambda}}{(n-1)!n^{k}}x^{n}, (\lambda \in \mathbb{R}),
\end{split}
\end{equation}
where $(x)_{0,\lambda} = 1$, $(x)_{n,\lambda} = x(x-\lambda)\cdots(x-(n-1)\lambda)$, $(n \ge 1)$ .\\
From \eqref{04}, we note that
\begin{equation}\label{05}
\begin{split}
Ei_{1,\lambda}(x) = \sum_{n=1}^{\infty}\frac{(1)_{n,\lambda}}{n!}x^{n} = e_{\lambda}(x)-1.
\end{split}
\end{equation}

The degenerate exponential functions are given by $e^{x}_{\lambda}(t) = (1+\lambda t)^{\frac{x}{\lambda}}$, $e_{\lambda}(t) = e_{\lambda}^{1}(t)$, and the degenerate logarithm functions are defined by $\log_{\lambda}(t)=\frac{1}{\lambda}(t^{\lambda}-1)$, which is the compositional inverse of $e_{\lambda}(t)$.\\

The degenerate poly--Genocchi polynomials are defined by
\begin{equation}\begin{split}\label{06}
\frac{2Ei_{k,\lambda}(\log_{\lambda}(1+t))}{e_{\lambda}(t)+1} e_{\lambda}^{x}(t) = \sum_{n=0}^{\infty} g_{n,\lambda}^{(k)}(x)\frac{t^{n}}{n!}, (k \in \mathbb{Z}).
\end{split}\end{equation}
For $x=0$, $g_{n,\lambda}^{(k)} =  g_{n,\lambda}^{(k)}(0)$ are called the degenerate poly--Genochhi numbers.\\
Here $g_{n,\lambda}^{(1)}(x)  = G_{n,\lambda}(x)$ are the degenerate Genocchi polynomials given by
\begin{equation}\label{07}
\begin{split}
\frac{2t}{e_{\lambda}(t)+1} e^{x}_{\lambda}(t) = \sum_{n=0}^{\infty}G_{n,\lambda}(x)\frac{t^{n}}{n!}, \quad (\textnormal{see}\,\,[12,13]).
\end{split}
\end{equation}
More generally, for $r \in \mathbb{N}$, the Genocchi polynomials of order $r$ are defined by
\begin{equation}\label{08}
\begin{split}
\left(\frac{2t}{e_{\lambda}(t)+1}\right)^{r}e_{\lambda}^{x}(t) = \sum_{n=0}^{\infty}G_{n,\lambda}^{(r)}(x)\frac{t^{n}}{n!}, \quad
(\textnormal{see}\,\,[13]).
\end{split}
\end{equation}
From \eqref{06}, we note that
\begin{equation}\label{09}
\begin{split}
g_{n,\lambda}^{(k)}(x) = \sum_{l=0}^{n}{n \choose l}g_{l,\lambda}^{(k)}(x)_{n-l,\lambda}, \quad (n \ge 0).
\end{split}
\end{equation}

We will need the Carlitz's degenerate Euler polynomials $\mathcal{E}^{(r)}_{l,\lambda}(x)$ of order $r$ given by
\begin{equation}\label{10}
\begin{split}
\left(\frac{2}{e_{\lambda}(t)+1}\right)^{r}e_{\lambda}^{x}(t) = \sum_{n=0}^{\infty} \mathcal{E}^{(r)}_{n,\lambda}(x)\frac{t^n}{n!}.
\end{split}
\end{equation}

For $k_{1}, k_{2},\cdots, k_{r} \in \mathbb{Z}$, we define the degenerate multiple polyexponential function as
\begin{equation}\label{10-1}
\begin{split}
Ei_{k_{1}, k_{2},\cdots, k_{r},\lambda} (x) = \sum_{0 <n_{1}<n_{2}<\cdots <n_{r}}\frac{(1)_{n_{1},\lambda\cdots}(1)_{n_{r},\lambda}x^{n_{r}}}{(n_{1}-1)!(n_{2}-1)!\cdots(n_{r}-1)!n_{1}^{k_{1}}\cdots n_{r}^{k_{r}}},
\end{split}
\end{equation}
where the sum is over all integers $n_1,n_2, \dots n_r$, satisfying $0 <n_{1}<n_{2}<\cdots <n_{r}$.

\medskip

\section{Degenerate multi-poly-Genocchi polynomials}
\medskip
For $\lambda \in \mathbb{R}$, the degenerate Stirling numbers of the first kind are defined by
\begin{equation}\label{11}
\begin{split}
(x)_{n} = \sum_{l=0}^{n}S_{1,\lambda}(n,l)(x)_{l,\lambda}, (n \ge 0), \quad (\textnormal{see}\,\,[1,2,5-18]),
\end{split}
\end{equation}
where $(x)_{0}=1$, $(x)_{n} = x(x-1)\cdots (x-n+1),\,\, (n \geq 1)$.\\
From \eqref{11}, we note that
\begin{equation}\label{12}
\begin{split}
\frac{1}{k!}\left(\log_{\lambda}(1+t)\right)^{k} = \sum_{n=k}^{\infty}S_{1,\lambda}(n,k)\frac{x^{n}}{n!}, \,\,(k \ge 0),
\end{split}
\end{equation}
and that $\lim_{\lambda \to 0}S_{1,\lambda} = S_{1}(n,l)$, where $S_{1}(n,l)$ is the Stirling number of the first kind.

For $k_{1}, k_{2},\cdots, k_{r} \in \mathbb{Z}$, we  consider the degenerate mulit--poly--Genocchi polynomials which are given by
\begin{equation}\label{13}
\begin{split}
\frac{2^{r}Ei_{k_{1}, k_{2},\cdots, k_{r,\lambda}(\log_{\lambda}(1+t))}}{(e_{\lambda}(t)+1)^{r}} e_{\lambda}^{x}(t) = \sum_{n=0}^{\infty}g_{n,\lambda}^{(k_{1}, k_{2},\cdots, k_{r})}(x)\frac{t^{n}}{n!}.
\end{split}
\end{equation}
For $x=0$, $g_{n,\lambda}^{(k_{1}, k_{2},\cdots, k_{r})} = g_{n,\lambda}^{(k_{1}, k_{2},\cdots, k_{r})}(0)$ are called the degenerate multi--poly--Genocchi numbers.\\
From \eqref{13}, we note that
\begin{equation}\label{15}
\begin{split}
g_{n,\lambda}^{(k_{1},\cdots, k_{r})}(x) = \sum_{l=0}^{n}{n \choose l}g_{l,\lambda}^{(k_{1}, k_{2},\cdots, k_{r})}(x)_{n-l,\lambda}, (n \ge 0).
\end{split}
\end{equation}

From \eqref{13}, we have
\begin{equation}\begin{split}\label{16}
&\sum_{n=0}^{\infty}g_{n,\lambda}^{(k_{1}, k_{2},\cdots, k_{r})}(x)\frac{t^{n}}{n!} \\
& = \left(\frac{2}{e_{\lambda}(t)+1}\right)^{r}e_{\lambda}^{x}(t)\sum_{0 <n_{1}<n_{2}<\cdots <n_{r}}\frac{(1)_{n_{1},\lambda\cdots}(1)_{n_{r},\lambda}(\log_{\lambda}(1+t))^{n_{r}}}{(n_{1}-1)!\cdots(n_{r}-1)!n_{1}^{k_{1}}\cdots n_{r}^{k_{r}}} \\
& = \sum_{l=0}^{\infty}\mathcal{E}_{l,\lambda}^{(r)}(x)\frac{t^{l}}{l!} \\
& \quad \times \sum_{0 <n_{1}<n_{2}<\cdots <n_{r}}\frac{(1)_{n_{1},\lambda\cdots}(1)_{n_{r},\lambda}}{(n_{1}-1)!\cdots(n_{r-1}-1)!n_{1}^{k_{1}}\cdots n_{r-1}^{k_{r-1}}n_{r}^{k_{r}-1}}\frac{\left(\log_{\lambda}(1+t)\right)^{n_{r}}}{n_{r}!} \\
& = \sum_{l=0}^{\infty}\mathcal{E}_{l,\lambda}^{(r)}(x)\frac{t^{l}}{l!} \\
& \quad \times \sum_{0 <n_{1}<n_{2}<\cdots <n_{r}\le m}\frac{(1)_{n_{1},\lambda\cdots}(1)_{n_{r},\lambda}S_{1,\lambda}(m,n_{r})}{(n_{1}-1)!\cdots(n_{r-1}-1)!n_{1}^{k_{1}} ,\cdots n_{r-1}^{k_{r-1}},n_{r}^{k_{r}-1}}\frac{t^{m}}{m!} \\
& = \sum_{n=r}^{\infty}\left(\sum_{l=0}^{n-r}{n \choose l}\mathcal{E}^{(r)}_{l,\lambda}(x)\sum_{0 <n_{1}<n_{2}<\cdots <n_{r}\le n-l}\frac{(1)_{n_{1},\lambda\cdots}(1)_{n_{r},\lambda}S_{1,\lambda}(n-l,n_{r})}{(n_{1}-1)!\cdots(n_{r-1}-1)!n_{1}^{k_{1}}\cdots n_{r-1}^{k_{r-1}}n_{r}^{k_{r}-1}}\right) \\
& \quad \times \frac{t^{n}}{n!}.
\end{split}\end{equation}

Comparing the coefficients on both sides of \eqref{16}, we have the following theorem.

\begin{theorem}
For $k_{1},k_{2},\dots,k_{r}\in\mathbb{Z}$, and $n, r \in \mathbb{N}$ with $n \geq r$, we have
\begin{equation*}\begin{split}
&g_{n,\lambda}^{(k_{1}, k_{2},\cdots, k_{r})}(x) \\
& = \sum_{l=0}^{n-r}{n \choose l}\mathcal{E}^{(r)}_{l,\lambda}(x)\sum_{0 <n_{1}<n_{2}<\cdots <n_{r}\le n-l}\frac{(1)_{n_{1},\lambda\cdots}(1)_{n_{r},\lambda}S_{1,\lambda}(n-l,n_{r})}{(n_{1}-1)!\cdots(n_{r-1}-1)!n_{1}^{k_{1}}\cdots n_{r-1}^{k_{r-1}}n_{r}^{k_{r}-1}},
\end{split}\end{equation*}
Further, we have $g_{n,\lambda}^{(k_{1}, k_{2},\cdots, k_{r})}(x)=0$,  for $0 \le n \le r-1$.
\end{theorem}

It is immediate to show that
\begin{equation}
\begin{split}\label{18}
\bigg(\frac{2}{e_{\lambda}(t)+1} \bigg)^r e_{\lambda}^{x}(t)
& = \frac{1}{t^r} \bigg(\frac{2t}{e_{\lambda}(t)+1} \bigg)^r e_{\lambda}^{x}(t) \\
&= \frac{1}{t^r} \sum_{n=r}^{\infty} G_{n,\lambda}^{(r)}(x) \frac{t^{n}}{n!}
= \sum_{n=0}^{\infty} \frac{G_{n+r,\lambda}^{(r)}(x)} {r! \binom{n+r}{n}} \frac{t^{n}}{n!}. \\
\end{split}
\end{equation}
Thus, by \eqref{18}, we get
\begin{equation}\label{19}
\mathcal{E}_{n,\lambda}^{(r)} (x) =  \frac{1} {r! \binom{n+r}{n}} G_{n+r,\lambda}^{(r)}(x), (n, r \geq 0).
\end{equation}

Therefore, by \eqref{19}, we obtain the following corollary.
\begin{corollary}
For $k_{1},k_{2},\dots,k_{r}\in\mathbb{Z}$, and $n, r \in \mathbb{N}$ with $n \geq r$, we have
\begin{equation*}\begin{split}
 g_{n,\lambda}^{(k_1,\cdots,k_r)} (x)& = \sum_{l=0}^{n-r} \frac{\binom{n}{l}} {r! \binom{l+r}{l}}  G_{l+r,\lambda}^{(r)} (x) \\
&\quad \times \sum_{0 < n_1 < n_2 < \cdots < n_r \leq n-l} \frac{(1)_{n_1,\lambda} \cdots (1)_{n_r,\lambda} S_{1,\lambda}(n-l,n_r)} {(n_1-1)! \cdots (n_{r-1}-1)! n_{1}^{k_1} \cdots n_{r-1}^{k_{r-1}} n_{r}^{k_{r}-1}}.
\end{split}\end{equation*}
\end{corollary}

From \eqref{13}, we note that
\begin{equation}
\begin{split}\label{21}
& \sum_{n=0}^{\infty} g_{n,\lambda}^{(k_1,\cdots,k_r)} (r) \frac{t^{n}}{n!}
= 2^r \bigg(1 - \frac{1}{2}\frac{2}{e_{\lambda}(t)+1} \bigg)^r Ei_{k_1,\cdots,k_r,\lambda}(\log_{\lambda} (1+t)) \\
& = \sum_{m=0}^{\infty} \bigg(\sum_{l=0}^{r} \binom{r}{l} (-1)^l 2^{r-l} \mathcal{E}_{m,\lambda}^{(l)} \bigg) \frac{t^{m}}{m!} \\
&\qquad \times \sum_{0 < n_1 < n_2 < \cdots < n_r} \frac{(1)_{n_1,\lambda} \cdots (1)_{n_r,\lambda} (\log_{\lambda} (1+t))^{n_r} } {(n_1-1)! \cdots (n_{r-1}-1)! n_{1}^{k_1} \cdots n_{r-1}^{k_{r-1}} n_{r}^{k_{r}-1} n_r!}
\end{split}
\end{equation}
\begin{equation*}
\begin{split}
& = \sum_{m=0}^{\infty} \bigg(\sum_{l=0}^{r} \binom{r}{l} (-1)^l 2^{r-l} \mathcal{E}_{m,\lambda}^{(l)} \bigg) \frac{t^{m}}{m!} \\
&\qquad \times \sum_{0 < n_1 < n_2 < \cdots < n_r \le j} \frac{(1)_{n_1,\lambda} \cdots (1)_{n_r,\lambda} S_{1,\lambda}(j,n_r)} {(n_1-1)! \cdots (n_{r-1}-1)! n_{1}^{k_1} \cdots n_{r-1}^{k_{r-1}} n_{r}^{k_{r}-1}} \frac{t^{j}}{j!} \\
& = \sum_{n=r}^{\infty} \bigg( \sum_{m=0}^{n-r} \sum_{l=0}^{r} \sum_{0 < n_1 < n_2 < \cdots < n_r \leq n-m} \\
&\qquad \frac{(1)_{n_1,\lambda} \cdots (1)_{n_r,\lambda} \binom{r}{l} \binom{n}{m} (-1)^l 2^{r-l} \mathcal{E}_{m,\lambda}^{(l)} S_{1,\lambda}(n-m,n_r)} {(n_1-1)! \cdots (n_{r-1}-1)! n_{1}^{k_1} \cdots n_{r-1}^{k_{r-1}} n_{r}^{k_{r}-1}} \bigg) \frac{t^{n}}{n!}.
\end{split}
\end{equation*}

Therefore, by \eqref{21}, we obtain the following theorem.

\begin{theorem}
For $k_{1},k_{2},\dots,k_{r}\in\mathbb{Z}$, and $n, r \in \mathbb{N}$ with $n \geq r$, we have
\begin{equation*}\begin{split}
& g_{n,\lambda}^{(k_1,\cdots,k_r)} (r) \\
& =  \sum_{m=0}^{n-r} \sum_{l=0}^{r} \sum_{0 < n_1 < n_2 < \cdots < n_r \leq n-m} \frac{(1)_{n_1,\lambda} \cdots (1)_{n_r,\lambda} \binom{r}{l} \binom{n}{m} (-1)^l 2^{r-l} \mathcal{E}_{m,\lambda}^{(l)} S_{1,\lambda}(n-m,n_r)} {(n_1-1)! \cdots (n_{r-1}-1)! n_{1}^{k_1} \cdots n_{r-1}^{k_{r-1}} n_{r}^{k_{r}-1}}.
\end{split}\end{equation*}
\end{theorem}

By \eqref{13}, we get
\begin{equation}
\begin{split}\label{23}
\sum_{n=0}^{\infty} g_{n,\lambda}^{(k_1,\cdots,k_r)} (x+y) \frac{t^{n}}{n!}
&= \frac{2^r Ei_{k_1,\cdots,k_r,\lambda}(\log_{\lambda} (1+t)) } {(e_{\lambda}(t)+1)^r} e_{\lambda}^{x}(t) e_{\lambda}^{y}(t) \\
& = \sum_{l=0}^{\infty} g_{l,\lambda}^{(k_1,\cdots,k_r)} (x) \frac{t^{l}}{l!} \sum_{m=0}^{\infty} (y)_{m,\lambda} \frac{t^{m}}{m!} \\
& = \sum_{n=0}^{\infty} \bigg(\sum_{l=0}^{n} \binom{n}{l} g_{l,\lambda}^{(k_1,\cdots,k_r)} (x) (y)_{n-l,\lambda} \bigg) \frac{t^{n}}{n!}.
\end{split}
\end{equation}

Thus, by comparing the coefficients on both sides of \eqref{23}, we get the following proposition.
\begin{proposition}
For $k_{1},k_{2},\dots,k_{r}\in\mathbb{Z}$, and any nonnegative integer $n$, we have
\begin{equation*}
g_{n,\lambda}^{(k_1,\cdots,k_r)} (x+y) = \sum_{l=0}^{n} \binom{n}{l} g_{l,\lambda}^{(k_1,\cdots,k_r)} (x) (y)_{n-l,\lambda}.
\end{equation*}
\end{proposition}

\section{Conclusion}
As we mentioned in the Introduction, studying various versions of some special polynomials and numbers has drawn attention of some mathematicians, and many interesting results about those polynomials and numbers have been obtained. To state a few, these include the degenerate Stirling numbers of the first and second kinds, degenerate central factorial numbers of the second kind, degenerate Bernoulli numbers of the second kind, degenerate Bernstein polynomials, degenerate Bell numbers and polynomials, degenerate central Bell numbers and polynomials, degenerate complete Bell polynomials and numbers, degenerate Cauchy numbers, and so on (see [11,12,14] and the references therein). We note here that the study has been carried out by using several different tools, such as generating functions, combinatorial methods, $p$-adic analysis, umbral calculus, differential equations, probability theory and so on. \\
\indent In this paper, we intorduced the degenerate multiple polyexponential functions which are multiple versions of the degenerate modified polyexponential functions. Then we defined the degenerate multi-poly-Gennocchi polynomials by means of those functions. In addition, we derived some explicit expressions for the degenerate multi-poly-Gennocchi polynomials and certain properties related to those polynomials. \\
\indent It is one of our future projects to continue this line of research, namely to study degenerate versions of certain special polynomials and numbers and to find their applications in physics, science and engineering as well as in mathematics.

\vspace{0.3in}

 {\bf {Availability of data and material:}} Not
applicable.

\vspace{0.1in}

{\bf {Competing interests:}} The authors declare no conflict of
interest.

\vspace{0.1in}

{\bf{Funding:}} Not applicable.

\vspace{0.1in}

{\bf{Author Contributions:}} D.S.K. and T.K. conceived of the
framework and structured the whole paper; D.S.K. and T.K. wrote
the paper; J.K. and H.Y.K. typed the paper; D.S.K. and
T.K.completed the revision of the article. All authors have read
and agreed to the published version of the manuscript.
\vspace{0.1in}

{\bf{Acknowledgements:}}The authors would like to thank Jangjeon
Institute for Mathematical Science for the support of this
research.

\medskip

\end{document}